\magnification=1200
\input amssym.def
\input epsf
\def \BI {{\Bbb I}}
\def \Sthree {{\Bbb S}^3}
\def\sqr#1#2{{\vcenter{\vbox{\hrule height.#2pt
     \hbox{\vrule width.#2pt height#1pt \kern#1pt
     \vrule width.#2pt} \hrule height.#2pt}}}}
\def\square{\ ${\mathchoice\sqr34\sqr34\sqr{2.1}3\sqr{1.5}3}$}

\centerline {\bf FOUR OBSERVATIONS ON $n$-TRIVIALITY AND BRUNNIAN LINKS}
\bigskip
\centerline {Theodore B. Stanford \footnote {} 
{Partial support for this paper was provided
by the Naval Academy Research Council.}}
\centerline {Mathematics Department}
\centerline {United States Naval Academy}
\centerline {572 Holloway Road}
\centerline {Annapolis, MD\ \ 21402}
\medskip
\centerline {\tt stanford@nadn.navy.mil}
\bigskip
\bigskip

\leftskip=1in
\rightskip=1in plus1fil
\noindent
{\bf Abstract.}\ \ 
Brunnian links have been known for a long time in knot theory, 
whereas the
idea of $n$-triviality is a recent innovation.  We illustrate
the relationship between the two concepts with 
four short theorems.
\bigskip
\bigskip

\leftskip=0in
\rightskip=0pt plus1fil
In 1892, Brunn introduced some nontrivial links with the
property that deleting any single component produces a
trivial link.  Such links are now called Brunnian
links. (See Rolfsen [7]).  Ohyama [5] introduced the
idea of a link which can be independantly undone in $n$
different ways.  Here ``undo'' means to change some set of
crossings to make the link trivial.  ``Independant'' means
that once you change the crossings in any one of the $n$
sets, the link remains trivial no matter what you do to the
other $n-1$ sets of crossings.  Philosophically, the ideas
are similar because, after all, once you delete one
component of a Brunnian link the result is trivial no matter
what you do to the other components.  We shall prove four
theorems that make the relationship between Brunnian links
and $n$-triviality more precise.

We shall show (Theorem~1) that an $n$-component Brunnian
link is $(n-1)$-trivial; (Theorem~2) that an $n$-component
Brunnian link with a homotopically trivial component is
$n$-trivial; (Theorem~3) that an $(n-k)$-component link
constructed from an $n$-component Brunnian link by twisting
along $k$ components is $(n-1)$-trivial; and (Theorem~4)
that a knot is $(n-1)$-trivial if and only if it is ``locally
$n$-Brunnian equivalent'' to the unknot.  At the end of the
paper we sketch a proof of Theorem~G, which generalizes
Theorems~1--3.

The property of $n$-triviality is closely related to
Vassiliev invariants. It is not hard to show that if a link
(or braid, string link, knotted graph, etc) is $n$-trivial,
then its Vassiliev invariants of order $<n$ vanish.  Also,
it follows from the work of a number of different authors
that a knot is $n$-trivial if and only if its Vassiliev
invariants of order $<n$ vanish.  (See for example 
[8] or Habiro [2].)

A {\it link} will be a tame, oriented link in oriented $\Sthree$.
We shall assume that the components of a link are ordered,
although the choice of ordering is essentially irrelevant.
Equivalence of links is up to the ambient
isotopy, and we work with regular diagrams in the usual
fashion.  A link of $n$ components is {\it Brunnian} if
every $(n-1)$-component sublink is trivial.  If $L$ is
an $n$-component link and $T \subset \{1,2,\dots n\}$,
then we denote by $L_T$ the sublink of $L$ obtained by
deleting the components with indices in $T$.  Thus 
an $n$-component link
$L$ is Brunnian if and only if $L_T$ is trivial for all
nonempty $T \subset \{1,2, \dots n\}$.

Suppose a link $L$ has a diagram with $n$ disjoint sets of
crossings $S_1, S_2, \dots S_n$.  (The $n$ here is
not necessarily related to the number of components.)
If $T \subset \{1,2, ... n\}$,
then we denote by $L(T)$ the link obtained from $L$ by
changing all the crossings in $\cup_{i \in T} S_i$.
The link is said to be {\it $n$-trivial} if it has
such a diagram with $L(T)$ trivial for all $T \ne \emptyset$.
Note that $n$-trivial implies $(n-1)$-trivial for
$n>0$.  (In some of the literature, $n$-trivial is
defined to be what we call here $(n+1)$-trivial.)
Figure~1 shows two $2$-trivial links, the
Borromean Rings and the Whitehead Link.
One possible way to choose the 
sets $S_i$ is indicated with
letters ``A'' and ``B''.

If $L$ is a link, then we denote the mirror image of
$L$ by $\hat L$.  If $T \subset \{1,2, \dots n\}$, then
we set $\overline T = \{1,2, \dots n\} - T$.

\bigskip 
\centerline {
\epsfxsize = 2 true in
\epsffile {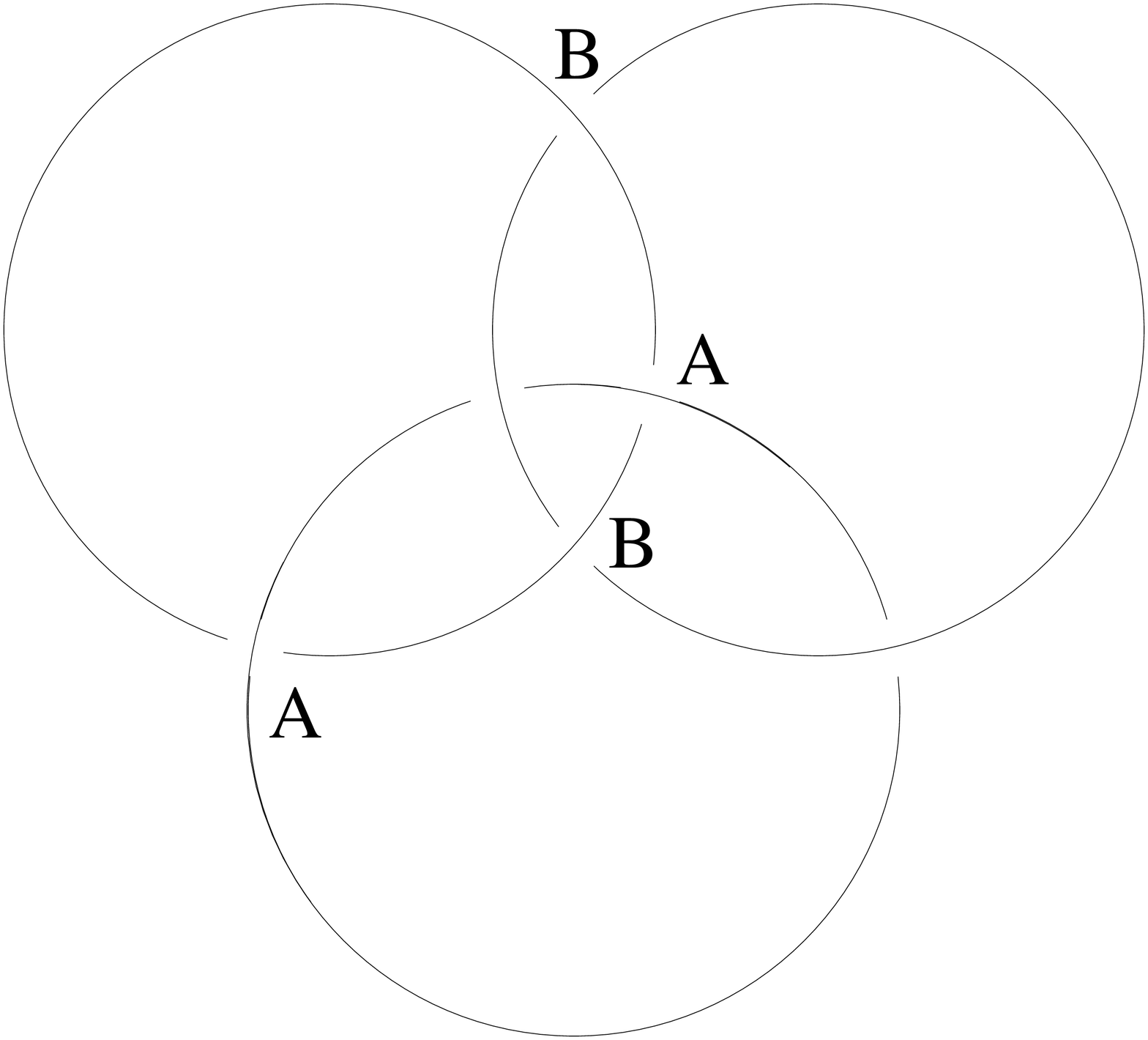}
\hskip 1 true in
\epsfxsize = 2 true in
\epsffile {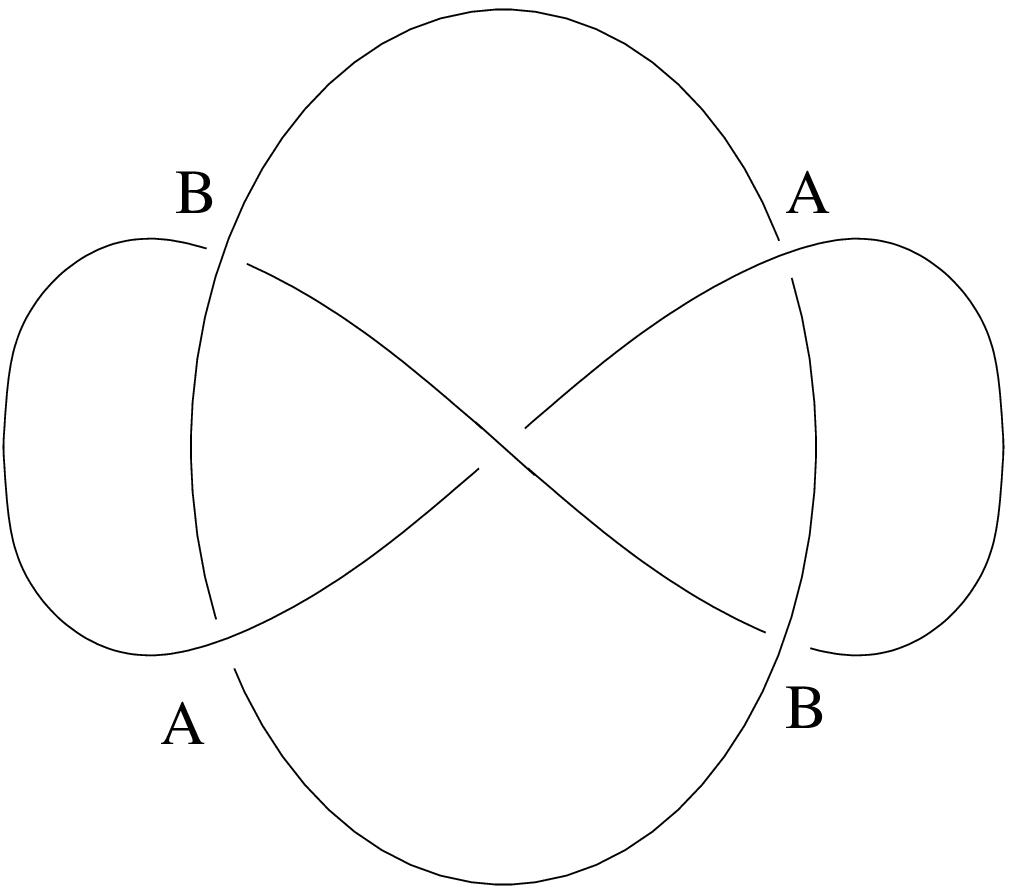}}
\medskip
\centerline {Borromean Rings \hskip 1.8 true in Whitehead Link}
\medskip
\centerline {Figure~1}
\bigskip

\medskip
\noindent
{\bf Theorem 1.}\ \ 
{\it An $n$-component Brunnian link is $(n-1)$-trivial.}
\medskip

{\it Proof:}\ \ 
Take any diagram of an $n$-component Brunnian link $L$.
For $1 \le i \le n$, let $S_i$ be the set of crossings where
the bottom strand is on the $i$th component.
(In the Borromean Rings in Figure~1, the crossings
labeled ``A'' and ``B'' correspond to $S_1$ and
$S_2$ chosen in this way.)
We need to show that if $T \subset \{1,2, \dots n-1\}$ and
$T \ne \emptyset$ then $L(T)$ is trivial.
In fact we will see that if $T$ is any proper
subset of $\{1,2, \dots n\}$, then $L(T)$ is trivial.
Let $T \subset \{1,2, \dots n\}$.
In $L(T)$, a strand from component $i \in T$ will always
pass over a strand from component $j \in \overline T$.
Thus $L(T)$ is the disjoint union of an
``upper'' link and a ``lower'' link.
The lower link is just $L_T$, since no crossings between
components with indices in $\overline T$ are changed. 
$L_T$ is a trivial link because $T \ne \emptyset$.
The upper link is $\hat L_{\overline T}$, consisting
of all components with indices in $T$, with all crossings
changed.  Since $L_{\overline T}$ is trivial (because
$\overline T \ne \emptyset$), so is $\hat L_{\overline T}$,
and therefore so is $L(T)$.
\square
\medskip

In particular, an $n$-component Brunnian link $L$ has trivial
Vassiliev invariants of order $<n-1$.  This was shown
by Kalfagianni and Lin [3] in the special case where
$L$ is the plat-closure of a pure braid.

The proof of Theorem 1 suggests that an $n$-component
Brunnian link is almost $n$-trivial, the only failure of the
sets $S_1, S_2, \dots S_n$ being when $T = \{1,2, \dots
n\}$.  We cannot hope for an $n$-component Brunnian link to
be $n$-trivial in general---this can fail in fact when
$n=2$.  Consider a two-component link with each component
unknotted and a nonzero linking number between the two
components. Such a link is
Brunnian, but if it were $2$-trivial then it would have
vanishing Vassiliev invariants of order $<2$, and it is
well-known that the linking number between two components of
a link is a Vassiliev invariant of order $1$.

We can, however, add $S_n$ to the argument if we make an
additional assumption.  We shall say that a component $K$ of
a link is {\it homotopically trivial} if there is a homotopy
in the complement of $L - K$ taking $K$ to a trivial loop.
In terms of link diagrams, this means that $L$ may be
trivialized by Reidemeister moves and by allowing crossing
changes of $K$ with itself.
Of the two links in Figure~1, the components
of the Whitehead Link are both homotopically trivial,
whereas none of the components of the Borromean Rings are.

\medskip
\noindent
{\bf Theorem 2.}\ \ 
{\it If $L$ is an $n$-component Brunnian link
with a homotopically trivial component, then $L$ is
$n$-trivial.}
\medskip

\noindent
{\it Proof:}\ \ Suppose without loss of generality that the
$n$th component $K$ is homotopically trivial.  Choose a
diagram for $L$ such that there exists a set $R$ of
crossings between strands of $K$ such that changing all of
the crossings in $R$ trivializes $L$.  Let $S_1, S_2, \dots
S_{n-1}$ be as before.  Let $S_n$ be the set of all
crossings where the bottom strand is on the $n$th component, minus
the crossings in $R$. 
(In the Whitehead link in Figure~1, the crossings
labeled ``A'' and ``B'' correspond to $S_1$ and
$S_2$ chosen in this way.)
Let $T \subset \{1,2, \dots n\}$ be nonempty.
As before, $L(T)$ is now the disjoint union of
a lower link and an upper link.  The lower link is
again $L_T$, and is always trivial (since $T \ne \emptyset$)
but it may be empty.  In the upper link, all crossings
of $L_{\overline T}$ are changed except those in $R$.
Therefore, the upper link is $\hat L_{\overline T}$ with the
crossings in $R$ changed, which is trivial for 
$T = \{1,2, \dots n\}$ by the choice of $R$.  For any
other $T$, the upper link will be a sublink of this
trivial link, and will therefore be trivial.
\square
\medskip

Let $L$ be an $n$-component Brunnian link, and let $U$ be a
proper subset of $\{1,2, \dots n\}$.  Since $L_{\overline U}$ is
trivial, its components bound disjoint disks.  Choose a
framing (an integer) for each component of $L_{\overline
U}$.  Let $L^U$ be obtained from $L_U$ by twisting along the
disk of each component of $L_{\overline U}$ according to its
framing.  Note that the components of $L^U$, as with $L_U$,
are the components of $L$ whose indices are in $\overline
U$.

\medskip
\noindent
{\bf Theorem 3.}\ \ {\it Let $L$ be an $n$-component
Brunnian link, and let $U \subset \{1,2, \dots n\}$.  Let
$L^U$ be obtained by twisting along the components of $L$
(using a fixed but arbitrary framing) whose indices are in
$U$, as above.  Then $L^U$ is $(n-1)$-trivial.}
\medskip

\noindent
{\it Proof:}\ \ 
We may assume without loss of generality
that $U = \{1,2, \dots k\}$ for some $k<n$.  Choose a diagram
for $L$ such that components $1,2, \dots n-1$ are all
disjoint circles.  Arrange the diagram, moreover, so that
around the $i$th component, $1 \le i \le k$, the diagram
looks like the left-hand side of Figure~2.  In general
there will be an arbitrary number of strands from the
$n$th component passing through, not just the three shown.
Now we may draw a diagram for $L^U$ by replacing each
local picture around the $i$th component, $1 \le i \le k$,
by the right-hand side of Figure~2.  A single twist is
shown, but there is a similar diagram for any integer
number of twists.  The point is that the twisted strands
can always be drawn such that there is a set of crossings,
like those marked with an ``X'' in Figure~2, such
that changing those crossings undoes the effect of the
twist.  For $1 \le i \le k$, let $S_i$ be that set of
crossings.  

For $k < i < n$, let $S_i$ be the set of crossings where the
lower strand is on the $i$th component, as before.  Because of the
way we have chosen the diagram, the top strands of these
crossings will all be on the $n$th component.  Now observe that
for $1 \le i < n$, changing the crossings in $S_i$ has the
same effect as removing the $i$th component from $L$ (for $k
< i < n$, ``remove'' means separate into a disjoint union)
and then twisting along any components in $U$ which may be
left.  But once one component has been removed from $L$ it
becomes trivial, and after twisting along or removing any
other components, it is still trivial.
\square
\medskip

\bigskip 
\centerline {
\epsfxsize = 5 true in
\epsfysize = 2 true in
\epsffile {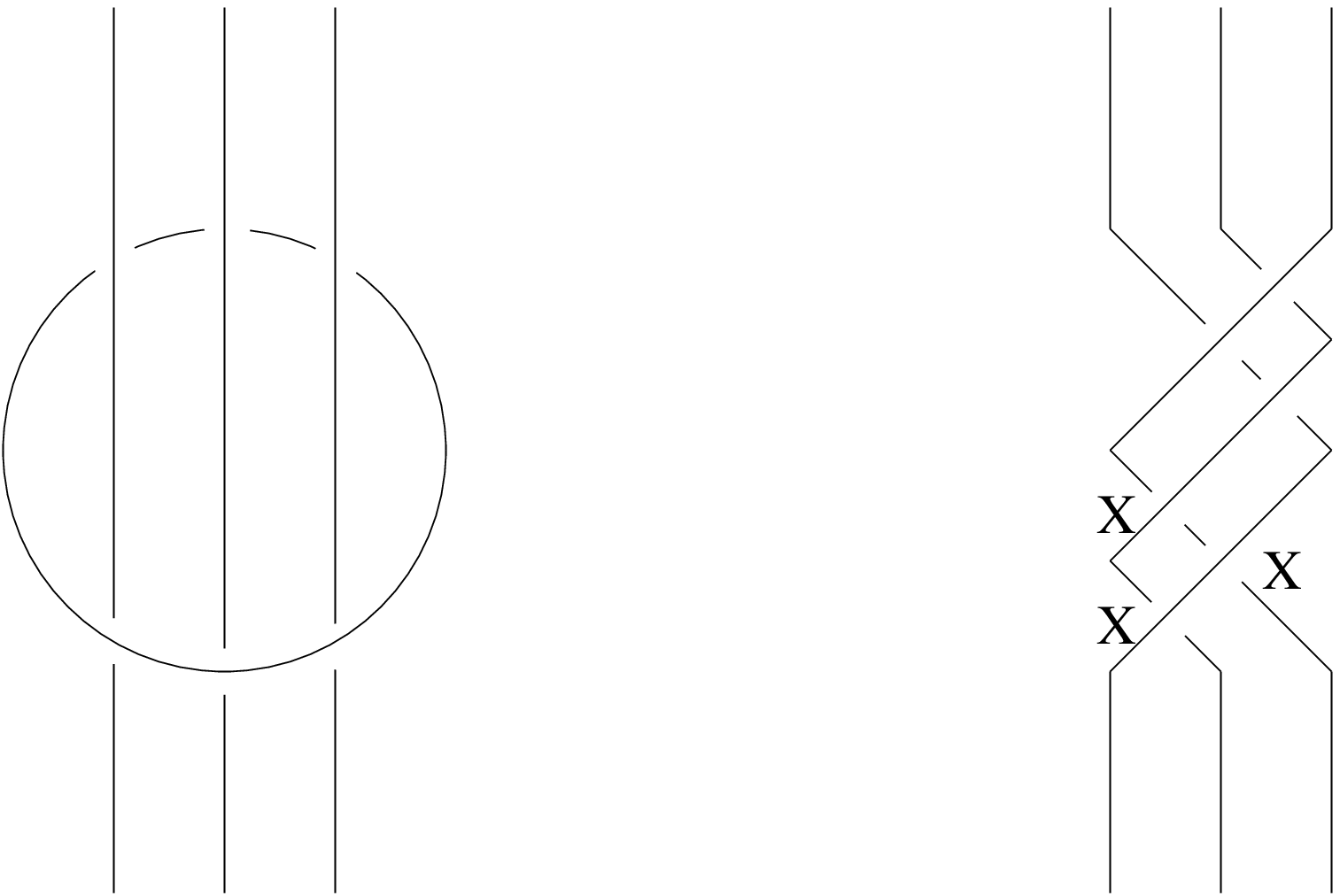}}
\medskip
\centerline {Figure~2}
\bigskip

As an example of Theorem~3,
it is not hard to see that twisting once along one component
of the Borromean Rings produces the Whitehead Link.
Both links are $2$-trivial, as noted above.

There have been several notions of $n$th-order equivalence
introduced for knots and links in the last few years.  There
are $n$-equivalence, $n$-similarity, $V_n$-equivalence, and
others.  All these notions are now known to be equivalent
for knots.  See Gusarov [1], Habiro [2], Ohyama [6],
Ng and Stanford [4], and [8].  We will add yet one more
characterization of this same idea in terms of Brunnian
string links.

For each positive integer $n$, fix $n$ distinct, ordered
points $x_1, x_2, \dots x_n$ in the two-dimensional disk
$\BI^2$.  An $n$-component {\it string link} is a proper,
tame embedding $f_1, f_2, \dots f_n$ of $n$ disjoint copies
of the unit interval $\BI$ into $\BI^3 = \BI^2 \times \BI$
such that $f_i (0) = (x_i,0)$ and $f_i(1) = (x_i,1)$ for all
$1 \le i \le n$.  Equivalence is up to ambient isotopy, fixing
the boundary.
We work with regular diagrams, as with
knots and links.  
An $n$-component string link $L$ is said
to be {\it Brunnian} if deleting any single component
trivializes $L$ (ie, makes it planar).  The string link shown
in Figure~3 is Brunnian.

\bigskip 
\centerline {
\epsfxsize = 3 true in
\epsfysize = 2.5 true in
\epsffile {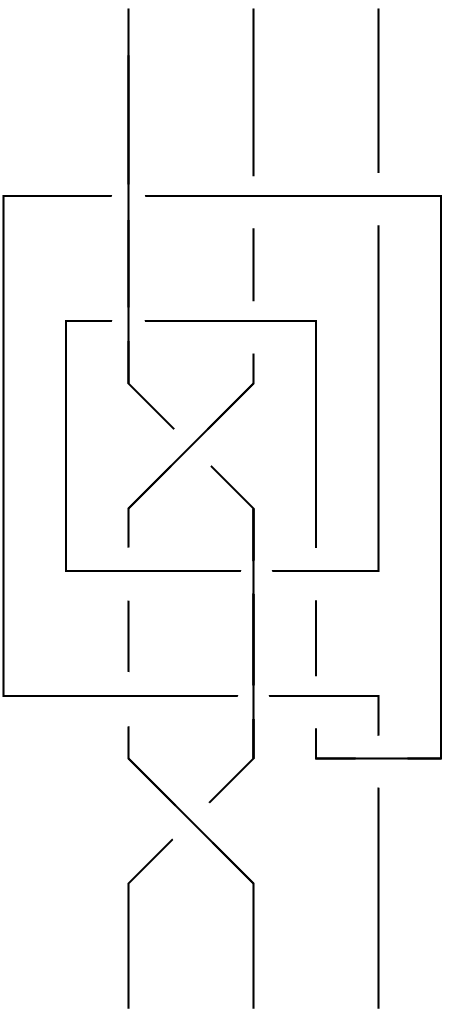}}
\medskip
\centerline {Figure~3}
\bigskip

\medskip
\noindent
{\bf Definition.}\ \ 
We say two knots are {\it locally $n$-Brunnian
equivalent} if one can be obtained from the other by a
sequence of local replacements of a trivial string link by
an $n$-component Brunnian string link.  We allow different
string links at each
each replacement.

\medskip
\noindent
{\bf Theorem 4.}\ \ 
{\it
Two knots $K$ and $K^\prime$
are locally $n$-Brunnian equivalent if and only if
they are $(n-1)$-equivalent.}
\medskip

\noindent {\it Proof:}\ \ 
An easy modification of the proof of Theorem~1 shows that
if $K^\prime$ is obtained from $K$ by replacing an
$n$-component trivial string link with a Brunnian
string link, then the two knots are $(n-1)$-similar in
Taniyama's sense (see Ohyama [5]).  
For the converse, observe that
the $C_n$ moves and the $*^n$ moves defined by Habiro [2]
are both examples of replacing a trivial string link
with $(n+1)$ components by a Brunnian string link.
\square
\medskip

\noindent
{\bf Remark:}\ \ 
Habiro's results indicate that Theorem~4 is
valid for string links, but
that for links in general the various 
notions of $n$th-order equivalence
diverge in ways that are not well-understood yet.
\medskip

We finish with a general theorem of which Theorems
1--3 are special cases.  
First, we define a {\it mixed link} to be a string
link with some circle components added, as in Figure~4.
We allow the cases of only string components (a string link)
or only circle components (a link in the usual sense).
Then we define {\it coloring} as a generalization of
ordering.  (Our coloring has nothing to do with counting
representations into finite groups.)  An {\it $n$-color
link} is a link together with a color (a number between $1$
and $n$) assigned to each component.  We require that all $n$
colors be used.  A color may have string components, or
circle components, or both.  An $n$-color link is {\it
Brunnian} if deleting the components of any single color
produces a trivial link.  It is easy to produce $n$-color
Brunnian links which are not Brunnian in the usual sense.
For example, take any $n$-component nontrivial Brunnian link
$L$ and form $L^\prime$ by replacing any component of $L$ by
two parallel unlinked
copies of itself.  Then $L^\prime$ is not
Brunnian in the usual sense.  It is, however, $n$-color
Brunnian if the two parallel components are given the same
color and each other component is given its own color.

\bigskip 
\centerline {
\epsfxsize = 2 true in
\epsfysize = 1.3 true in
\epsffile {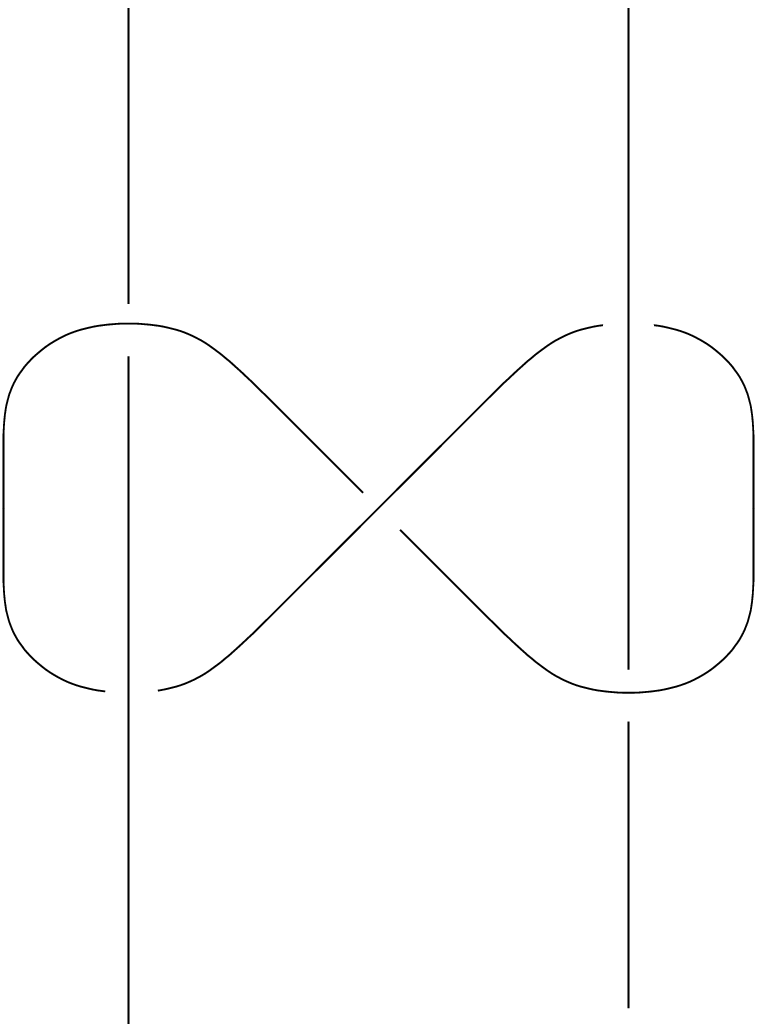}}
\medskip
\centerline {Figure~4}
\bigskip

Let $L$ be an $n$-color mixed 
link, and let $k \in \{1,2, \dots n\}$.
The color $k$ is said to be {\it homotopically trivial} if
there exists a homotopy of the components of color $k$ to
the trivial link in the complement of all the remaining
components.  In terms of diagrams, this means that $L$
can be trivialized by Reidemeister moves and crossing
changes between two strands both of color $k$.
For example, if the Borromean Rings in Figure~1 are
colored with two colors, then it is easy to see that 
the color with two components is homotopically trivial,
whereas no single component of the link is homotopically
trivial.

If $L$ is an $n$-color link and $T \subset \{1,2, \dots n\}$,
then we may define $L_T$ to be obtained from 
$L$ by deleting the components whose color is in $T$.
If $U \subset \{1,2, \dots n\}$, and if
all the components with color in $U$ are circle
components, then we may form $L^U$ as above by
choosing a framing for each component with color in $U$
and twisting along each such component according to its
chosen framing.  

\medskip
\noindent
{\bf Theorem G.}\ \ 
{\it 
Let $L$ be an $n$-color mixed link. Let $U$ be a 
subset of $\{1,2, \dots n\}$, possibly empty but
not equal to $\{1,2, \dots n\}$,
such that every component with color in $U$ is a circle
component.  Let $L^U$ be obtained from $L$ by twisting
along the components with colors in $U$ according to
some fixed but arbitrary set of framings.  Then $L^U$
is $(n-1)$-trivial.  Moreover, if one color of $L$ 
not in $U$ is homotopically trivial, then $L^U$ is
$n$-trivial.}
\medskip

The proof of Theorem~G is mostly a matter of putting together
the proofs of Theorems~1--3.  If there is a homotopically
trivial color, we may assume that it is $n$.  Choose a
diagram for $L$ which is planar on the sublink of the first
$n-1$ colors.  If the $n$th color is homotopically trivial,
then choose the diagram so that there exists a set $R$ of
crossings between strands of color $n$, such that changing
all the crossings in $R$ trivializes $L$.
Now draw the diagram for $L^U$ obtained by replacing a disk
around each component of color $i$, $1 \le i \le k$, with a
local twist (or with several twists) as in Figure~2.  
For $1 \le i \le k$, let $S_i$ be the set of
crossings in $L^U$ which undoes the effects of the twists
around all the components with color $i$.  For $k < i < n$, let
$S_i$ be the set of crossings whose bottom strand has color
$i$.  If the color $n$ is homotopically trivial, then let
$S_n$ be the set of crossings whose bottom strand has color
$n$, minus the crossings in $R$, minus any crossings
introduced by the twists along the first $k$ components.
Then one only has to check that all the appropriate links
$L^U(T)$ are trivial, as in the previous proofs.

\bigskip
\noindent
{\bf Acknowledgement.}\ \ 
I would like to thank John Dean for some helpful
and stimulating conversations.

\bigskip
\noindent
{\bf References.}
\medskip

\smallskip 
\item {1.} M. N. Gusarov.
On $n$-equivalence of knots and invariants of finite degree,
{\it Topology of Manifolds and Varieties}, 173--192, 
{\it Advances in Soviet Mathematics} 18, 
American Mathematical Society, 1994

\smallskip 
\item {2.} K. Habiro. 
Claspers and the Vassiliev skein modules.
Preprint, University of Tokyo.

\smallskip
\item {3.}
E. Kalfagianni and X.-S. Lin.
Milnor and finite type invariants of plat-closures.
Preprint {\tt GT/9804030} available from 
{\tt front.math.ucdavis.edu}.

\smallskip
\item {4.}
K.Y. Ng and T. Stanford.
On Gusarov's groups of knots.
To appear in {\it  Mathematical Proceedings of the 
Cambridge Philosophical Society}.

\smallskip
\item {5.} 
Y. Ohyama.
A new numerical invariant of knots induced from their
regular diagrams. {Topology and its Applications} 37 (1990)
no. 3, 249--255.

\smallskip
\item {6.} 
Y. Ohyama.
Vassiliev invariants and similarity of knots.
{Proceedings of the American Mathematical Society}
123 (1995) no. 1, 287--291.

\smallskip
\item {7.}
D. Rolfsen.
{\it Knots and Links}, volume~7 of 
{\it Mathematics Lecture Series}.
Publish or Perish, Inc., Wilmington, DE, 1976.

\smallskip
\item {8.}
T.B. Stanford.
Vassiliev invariants and knots modulo pure braid subgroups.
Preprint {\tt GT/9805092} available from 
{\tt front.math.ucdavis.edu}.

\end